\documentclass{amsart}
\usepackage{amssymb, amsmath, latexsym}

\renewcommand{\baselinestretch}{\baselinestretch}
\renewcommand{\baselinestretch}{1.1}
\author{Constantin-Nicolae Beli}
\title[]{Explicit formulas for the cohomology of elementary abelian
$p$-groups} 
\date{}
     
 \def\e{\varepsilon}  
   
   \def\m{\lim}

 \def\({\overline}

\def\){\underline} \def\<{\cdot} 
\def\>{~~~~~~~} \def\#{{\bf
Definition}} \def\*{\section} \def\be{\begin{equation}}
\def\ee{\end{equation}}

   \def\FF{{\mathbb F}}

\def\bmat{\left(\begin{array}} \def\emat{\end{array}\right)}

 \def\m2{~(\mo 2)} \def\no{\noindent}
 \def\btm{\begin{thm}}
\def\etm{\end{tm}}
 \def\blem{\begin{lem}}
\def\elem{\end{lem}}

\newtheorem{theorem}{Theorem}[section]
\newtheorem{proposition}[theorem]{Proposition}
\newtheorem{lemma}[theorem]{Lemma}
\newtheorem{definition}{Definition}
\newtheorem{corollary}[theorem]{Corollary}

\newtheorem{bof}[theorem]{}
\newtheorem{teorema}{Theorem}

\def\qed{\mbox{$\Box$}\vspace{\baselineskip}}
\def\pf{$Proof.$ } 
\def\bco{\begin{corollary}} \def\eco{\end{corollary}} 
\def\bdf{\begin{definition}} \def\edf{\end{definition}} 
\def\btm{\begin{theorem}} \def\etm{\end{theorem}} 
\def\bpr{\begin{proposition}} \def\epr{\end{proposition}}  
\def\blm{\begin{lemma}} \def\elm{\end{lemma}} 
\def\bff{\begin{bof}\rm} \def\eff{\end{bof}}
\def\btr{\begin{teorema}} \def\etr{\end{teorema}}
\def\II{{\mathcal I}}

\def\de{\newcommand} \de\tm[1]{{\no\bf Theorem~#1}} 
\def\mb{\mathbb} 
   \def\ZZ{{\mb Z}}
 \def\FF{{\mb F}}
\de\lm[1]{{\no\bf Lemma~#1}}
\de\df[1]{{\no\bf Definition~#1}} \de\co[1]{{\no\bf Corollary~#1}}
\de\lr[1]{\longrightarrow^{\!\!\!\!\!\!\!\! #1}}
\de\lf[1]{\longleftarrow^{\!\!\!\!\!\!\!\! #1}}
\de\si[1]{\sim^{\!\!\!\!\! #1}} \de\apr[1]{\approx^{\!\!\!\!\! #1}}
\de\leg[2]{\left(\frac {#1}{#2}\right)}
\DeclareMathOperator\Br{Br}

\de\Brr[1]{{}_{#1}\Br}

\DeclareMathOperator\sgn{sgn}
\DeclareMathOperator\Ima{Im}

\begin{document}

\begin{abstract}
Let $G$ be an elementary $p$-group, $G\cong\FF_p^r$, with the basis
$s_1,\ldots,s_r$, and let $V$ be its dual, $V={\rm
Hom}(G,\FF_p)=H^1(G,\FF_p)$. Let $x_1,\ldots,x_r$ be the basis of $V$
dual to $s_1,\ldots,s_r$ and let $y_i=\beta (x_i)\in H^2(G,\FF_p)$,
where $\beta :H^1(G,\FF_p)\to H^2(G,\FF_p)$ be the connecting
Bockstein map. Then the ring $(H^*(F,\FF_p),+,\cup )$ satisfies
$$H^*(G,\FF_p)\cong\begin{cases}\FF_2[x_1,\ldots,x_r]&p=2\\
\Lambda (x_1,\ldots,x_r)\otimes\FF_p[y_1,\ldots,y_r]&p>2\end{cases}.$$
If $p=2$ the isomorphism $\tau :\FF_2[x_1,\ldots,x_r]\to H^*(G,\FF_2)$ is
given by
$$x_{i_1}\cdots x_{i_n}\mapsto x_{i_1}\cup\cdots\cup
x_{i_n}\in H^n(G,\FF_2).$$
If $p>2$ the isomorphism $\tau :\Lambda 
(x_1,\ldots,x_r)\otimes\FF_p[y_1,\ldots,y_r]\to H^*(G,\FF_p)$ is given by
$$x_{i_1}\wedge\cdots\wedge x_{i_k}\otimes y_{j_1}\cdots y_{j_l}\mapsto
x_{i_1}\cup\cdots\cup x_{i_k}\cup y_{j_1}\cup\cdots\cup y_{j_l}\in
H^{k+2l}(G,\FF_p).$$

In this paper we give an exact formula for the reverse isomorphisms
$\tau^{-1}$ in the cases $p=2$ and $p>2$. The elements of
$H^*(G,\FF_p)$ are written in terms of normalized cochains.

The tool we use in our proof is an alernative desription of the
normalized cochains $C^n(G,M)\cong{\rm Hom}(T^n(\II ),M)$, where $\II
=\II_G$ is the augmentation ideal, $\II =\ker (\varepsilon :\ZZ
[G]\to\ZZ )$. This allows to write the maps $d_n:C^n(G,M)\to
C^{n+1}(G,M)$ in a very convenient form.

\end{abstract}
\maketitle

\section{The $\II$-cochains}

Let $G$ be a group. We denote by $\varepsilon :\ZZ[G]\to\ZZ$ the
augmentation map, given by $\sum_{s\in G}m_ss\mapsto\sum_{s\in
G}m_s$. Here $m_s$ are integers and (if $G$ is infinite) almost all
are zero. We denote by $\II =\II_G$ the augmented ideal,
$\II=\ker\e$. Then $s-1$, with $s\in G\setminus\{ 1\}$ are a basis of
$\II$ over $\ZZ$. 

If $M$ is a $G$-module, then for every $n\geq 0$ we denote by
$C^n(G,M)$ the normalized cochains of degree $n$, i.e.
$$C^n(G,M)=\{ a\in M^{G^n}\,\mid\, 
a(u_1,\ldots,u_n)=0\text{ if }u_i=1\text{ for some }i\}.$$

\bdf For every $n\geq 0$ we denote by $C_\II^n(G,M)={\rm Hom}\, (T^n(\II
),M)$. The elements of $C_\II^n(G,M)$ are called $\II$-cochains.
\edf

\blm For any $n\geq 0$ there is an one to one correspondence between
$C_\II^n(G,M)$ and $C^n(G,M)$, given by $f\mapsto a$, where
$$a(u_1,\ldots,u_n)=f((u_1-1)\otimes\cdots\otimes (u_n-1)).$$
\elm
\pf The elements $u-1$, $u\in G\setminus\{ 0\}$ are a basis of $\II$
over $\ZZ$ so $(u_1-1)\otimes\cdots\otimes (u_n-1)$, with $u_i\in
G\setminus\{ 1\}$ are a a basis for $T^n(\II )$. Therefore we have a
linear bijection $C_\II^n(G,M)={\rm Hom}(T^n(\II ),M)\to
M^{(G\setminus\{ 1\} )^n}$ given by $f\to\tilde a$, where $\tilde
a(u_1,\ldots,u_n)=f((u_1-1)\otimes\cdots\otimes (u_n-1)$ for
$(u_1,\ldots,u_n)\in (G\setminus\{ 1\} )^n$. Next, we have a bijection
$M^{(G\setminus\{ 1\} )^n}\to C^n(G,M)$, given by $\tilde a\mapsto a$,
where $a(u_1,\ldots,u_n)=\tilde a(u_1,\ldots,u_n)$ if 
$u_1,\ldots,u_n\in G\setminus\{ 1\}$ and $a(u_1,\ldots,u_n)=0$ when
$u_i=1$ for some $i$. By composing these bijections, we get the
claimed isomorphism $C_\II^n(G,M)\to C^n(G,M)$. The relation
$a(u_1,\ldots,u_n)=f((u_1-1)\otimes\cdots\otimes (u_n-1))$, which holds
when $(u_1,\ldots,u_n)\in (G\setminus\{ 1\} )^n$, remains true if
$u_i=1$, so $u_i-1=0$, for some $i$. In this case we have 
$f((u_1-1)\otimes\cdots\otimes
(u_n-1))=f(0)=0=a(u_1,\ldots,u_n)$. \qed

As a consequence of Lemma 1.1, $C_\II^n(G,M)$, instead of $C^n(G,M)$,
can be used as cochains for computing $H^*(G,M)$. It turns out that
the coboundary maps $d_n:C_\II^n(G,M)\to C_\II^{n+1}(G,M)$ can be
written in a very convenient form.

\bpr The coboundary map $d_n:C_\II^n(G,M)\to C_\II^{n+1}(G,M)$ is given by
$f\mapsto g$, where, for every $\alpha_1,\ldots\alpha_{n+1}\in\II$,
$$g(\alpha_1\otimes\cdots\otimes\alpha_{n+1})
=\alpha_1f(\alpha_2\otimes\cdots\otimes\alpha_{n+1})
+\sum_{i=1}^n(-1)^if(\alpha_1\otimes\cdots\otimes
\alpha_i\alpha_{i+1}\otimes\cdots\otimes \alpha_{n+1}).$$
\epr
\pf The $\II$-cochains $f\in C_\II^n(G,M)$ and $g\in C_\II^{n+1}(G,M)$
correspond to $a\in C^n(G,M)$ and $b\in C^{n+1}(G,M)$, where
$a(u_1,\ldots,u_n)=f((u_1-1)\otimes\cdots\otimes (u_n-1))$ and
$b(u_1,\ldots,u_{n+1})=g((u_1-1)\otimes\cdots\otimes (u_{n+1}-1))$
$\forall u_1,\ldots,u_{n+1}\in G$.

The condition $g=d_nf$, i.e. $b=d_na$, writes as
\begin{multline*}
b(u_1,\ldots,u_{n+1})=u_1a(u_2,\ldots,u_{n+1})+
\sum_{j=1}^n(-1)^ja(u_1,\ldots,u_ju_{j+1},\ldots,u_{n+1})\\
+(-1)^{n+1}a(u_1,\ldots,u_n).
\end{multline*}
We write $\alpha_i=u_i-1$ for $i=1,\ldots,n+1$. Then
$b(u_1,\ldots,u_{n+1})=g(\alpha_1\otimes\cdots\otimes\alpha_{n+1})$.

We now write the terms of the right hand side of the equation above in
terms of $f$. We have
$u_1a(u_2,\ldots,u_{n+1})=(\alpha_1+1)a(u_2,\ldots,u_{n+1})
=(\alpha_1+1)f(\beta_0)$, where
$$\beta_0=\alpha_2\otimes\cdots\otimes\alpha_{n+1}.$$
For $1\leq j\leq n$ we have
$(-1)^ja(u_1,\ldots,u_ju_{j+1},\ldots,u_{n+1})=f(\beta_j)$, where
\begin{align*}\beta_j & =(-1)^j(u_1-1)\otimes\cdots\otimes
(u_{j-1}-1)\otimes (u_ju_{j+1}-1)\otimes
(u_{j+2}-1)\otimes\cdots\otimes (u_{n+1}-1)\\
 & =(-1)^j\alpha_1\otimes\cdots\otimes
\alpha_{j-1}\otimes
(\alpha_j+\alpha_{j+1}+\alpha_j\alpha_{j+1})\otimes\alpha_{j+2}
\otimes\cdots\otimes\alpha_{n+1}.
\end{align*}
(We have $u_i-1=\alpha_i$ and $u_ju_{j+1}-1=(\alpha_j+1)(\alpha_{j+1}+1)-1
=\alpha_j+\alpha_{j+1}+\alpha_j\alpha_{j+1}$.)

Finally, $(-1)^{n+1}a(u_1,\ldots,u_n)=f(\beta_{n+1})$, where
$$\beta_{n+1}=(-1)^{n+1}\alpha_1\otimes\cdots\otimes\alpha_n.$$
By adding these relations, we get
$$g(\alpha_1\otimes\cdots\otimes\alpha_{n+1})=(\alpha_1+1)f(\beta_0)
+\sum_{j=1}^{n+1}f(\beta_j)=\alpha_1f(\beta_0)+f(\beta ),$$
where $\beta=\sum_{j=0}^{n+1}\beta_j$. Now, in the sum $\beta$ there are
terms of the form
$\alpha_1\otimes\cdots\otimes\hat\alpha_i\otimes\cdots\otimes\alpha_{n+1}$
with $1\leq i\leq n+1$. Any such term appears twice, once in
$\beta_{i-1}$ and once in $\beta_i$, with the coefficients
$(-1)^{i-1}$ and $(-1)^i$, which cancel each other. We are left with
terms like $\alpha_1\otimes\cdots\otimes
\alpha_i\alpha_{i+1}\otimes\cdots\otimes\alpha_{n+1}$, which appears
in $\beta_i$, with the coefficient $(-1)^i$. In conclusion, $\beta
=\sum_{i=1}^n(-1)^i\alpha_1\otimes\cdots\otimes
\alpha_i\alpha_{i+1}\otimes\cdots\otimes\alpha_{n+1}$. It follows that
$$g(\alpha_1\otimes\cdots\otimes\alpha_{n+1})=
\alpha_1f(\alpha_2\otimes\cdots\otimes\alpha_{n+1})+
f\Big(\sum_{i=1}^n(-1)^i\alpha_1\otimes\cdots\otimes
\alpha_i\alpha_{i+1}\otimes\cdots\otimes\alpha_{n+1}\Big).$$
Hence $d_nf=g$ is equivalent to the desired equality for
$\alpha_1,\ldots,\alpha_{n+1}$ of the form $u-1$, with $u\in G$. But
both sides of the relation we want to prove are linear in each of the
variables $\alpha_1,\ldots,\alpha_{n+1}\in\II$ and $\{ u-1\mid\,
u\,\in G\}$ generate $\II$. Hence our relation holds for
$\alpha_1,\ldots,\alpha_{n+1}\in\II$ arbitrary. \qed

\bff{\bf Remark.} If $a\in M^G$, in particular, if the action of $G$
on $M$ is trivial, then $(s-1)a=0$ $\forall s\in G$, which implies
that $\alpha a=0$ $\forall\alpha\in\II$. 

Consequently, if the action of $G$ on $M$ is trivial then
$$d_nf(\alpha_1\otimes\cdots\otimes\alpha_{n+1})=
\sum_{i=1}^n(-1)^if(\alpha_1\otimes\cdots\otimes \alpha_i\alpha_{i+1}
\otimes\cdots\otimes\alpha_{n+1}).$$
\eff

Unfortunately, in the general case, for the cup products we don't have
a nice formula such as for the coboundary map. We have however the
following partial result.

\blm (i) If $M,N$ are $G$-modules and the action of $G$ on $N$ is
trivial, then for every $m,n\geq 0$ the cup product $\cup
:C_\II^m(G,M)\times C_\II^n(G,N)\to C_\II^{m+n}(G,M\otimes N)$ is
given by
$$(f\cup g)(\alpha_1\otimes\cdots\otimes\alpha_{m+n})=
(-1)^{mn}f(\alpha_1\otimes\cdots\otimes\alpha_m)\otimes
g(\alpha_{m+1}\otimes\cdots\otimes\alpha_{m+n}).$$

(ii) More generally, if $M_1,\ldots,M_r$ are $G$-modules and the
action of $G$ on $M_2,\ldots,M_r$ is trivial, then for every
$n_1,\ldots,n_r\geq 0$, with $n_1+\cdots +n_r=n$, the cup product
$C_\II^{n_1}(G,M_1)\times\cdots\times C_\II^{n_r}(G,M_r)\to
C_\II^n(G,M_1\otimes\cdots\otimes M_r)$ is given by
$$(f_1\cup\cdots\cup f_r)(\alpha_1\otimes\cdots\otimes\alpha_n)
=(-1)^{\frac{l(l-1)}2}\bigotimes_{i=1}^rf_i(\alpha_{n_1+\cdots
+n_{i-1}+1}\otimes\cdots\otimes\alpha_{n_1+\cdots +n_i}),$$
where $l=|\{ i\,\mid\, 1\leq i\leq r,\, n_i\text{ is odd}\}|$.
\elm
\pf (i) Let $a\in C^m(G,M)$ and $b\in C^n(G,N)$ be the normalized
cochains corresponding to $f$ and $g$.

The relation we want to prove writes as $(f\cup
g)(\alpha_1\otimes\cdots\otimes\alpha_{m+n})
=h(\alpha_1,\ldots,\alpha_{m+n})$, where $h:\II^{m+n}\to\FF_p$ is
given by $h(\alpha_1,\ldots,\alpha_{m+n})
=(-1)^{mn}f(\alpha_1\otimes\cdots\otimes\alpha_m)\otimes
g(\alpha_{m+1}\otimes\cdots\otimes\alpha_{m+n})$. Since both sides of
this relation are linear in each variable, it suffices to consider the
case when $\alpha_1,\ldots,\alpha_{m+n}$ belong to the set of
generators $\{ u-1\mid u\in G\}$ of $\II$. If we take
$\alpha_i=u_i-1$, with $u_i\in G$, then our relation writes in terms
of the normalized cochains $a$ and $b$ as
$$(a\cup b)(u_1,\ldots,u_{m+n})=(-1)^{mn}a(u_1,\ldots,u_m)\otimes
b(u_{m+1},\cdots u_{m+n}).$$
But this follows from the definition of the cup product and the fact
that the action of $G$ on $N$ is trivial, so $u_1\cdots
u_mb(u_{m+1},\cdots u_{m+n})=b(u_{m+1},\cdots u_{m+n})$.

(ii) If we apply (i) to $M=M_1\cup\cdots\cup M_i$ and $N=M_{i+1}$ we
see that the factor $(-1)^{(n_1+\cdots+n_i)n_{i+1}}$ appears. Hence,
by induction, we have
$$(f_1\cup\cdots\cup f_r)(\alpha_1\otimes\cdots\otimes\alpha_n)
=(-1)^S\bigotimes_{i=1}^rf_i(\alpha_{n_1+\cdots
+n_{i-1}+1}\otimes\cdots\otimes\alpha_{n_1+\cdots +n_i}),$$
where $S=\sum_{i=1}^{r-1}(n_1+\cdots +n_i)n_{i+1}=\sum_{1\leq i<j\leq
r}n_in_j$. But for $i<j$ the product $n_in_j$ is odd iff both $i,j$
belong to the set $\{ i\,\mid\, 1\leq i\leq r,\, n_i\text{ is odd}\}$,
of cardinal $l$. There are $\binom l2=\frac{l(l-1)}2$ such pairs
$(i,j)$. So the sum $S$ has $\frac{l(l-1)}2$ odd terms so
$S\equiv\frac{l(l-1)}2\pmod 2$, which concludes the proof. \qed

We denote by $Z_\II^n(G,M)=\ker
(C_\II^n(G,M)\xrightarrow{d_n}C_\II^{n+1}(G,M))$, $B_\II^n(G,M)=\Ima
(C_\II^{n-1}(G,M)\xrightarrow{d_{n-1}}C_\II^n(G,M))$ and
$H^n(G,M)=Z_\II^n(G,M)/B_\II^n(G,M)$. The isomorphism between
$C_\II^n(G,M)$ and $C^n(G,M)$ induces isomorphisms between $Z^n(G,M)$,
$B^n(G,M)$ and $H^n(G,M)$ and $Z_\II^n(G,M)$,
$B_\II^n(G,M)$ and $H_\II^n(G,M)$, respectively. 
\medskip

{\bf Remark} The $\II$-cochains bear some similarities to the cochains
resulting from the {\em Alternative Description of the Bar Resolution}
from [HS, chapter VI, 13(c)] they are not related. For convenience, we
denote by $\bar C^n(G,M)$ the normalized homogeneous cochains and by
$\bar C_\II^n(G,M)$ the cochains resulting from the alternative
description of the bar resolution, i.e. $\bar C_\II^n(G,M)={\rm
Hom}_G(\ZZ G\otimes T^n(\II ),M)$ If $a\in C^n(G,M)$ and $\bar a$, $f$
and $\bar f$ are its correspondents in $\bar C^n(G,M)$, $C_\II^n(G,M)$
and $\bar C_\II^n(G,M)$, then $a(u_1,\ldots,u_n)=\bar
a(1,u_1,u_1u_2,\ldots,u_1\cdots u_n)$ and $\tilde
a(u_0,\ldots,u_n)=\tilde f(u_0\otimes (u_1-u_0)\otimes\cdots\otimes
(u_n-u_{n-1}))$ $\forall u_0,\ldots,u_n\in G$. Hence for every
$u_1,\ldots,u_n\in G$ we have
$$\begin{aligned}
f((u_1-1)\otimes\cdots\otimes (u_n-1))&=a(u_0,\ldots,u_n)=\tilde
a(1,u_1,u_1u_2,\ldots,u_1\cdots u_n)\\
&=\tilde f(1\otimes (u_1-1)\otimes u_1(u_2-1)\otimes\cdots\otimes
u_1\cdots u_{n-1}(u_n-1)).
\end{aligned}$$ 
As one can see, the relation between $f$ and $\tilde f$ is not as
simple as $f(\eta )=\tilde f(1\otimes\eta )$ $\forall\eta\in T^n(\II
)$. Instead, it is somewhat similar to the relation between $a$ and
$\bar a$.

\section{Main result}

Let $G\cong\FF_p^r$ where $p$ is a prime and $r\geq 1$. Let
$s_1,\ldots,s_r$ be a basis of $G$ over $\FF_p$. On $G$ we use the
multiplicative notation. Let $V=G^*={\rm Hom}_{\FF_p}(G,\FF_p)$ and
let $x_1,\ldots,x_r\in V$ be the dual basis for the basis
$s_1,\ldots,s_r$ of $G$. We have $x_i(s_j)=\delta_{i,j}$.

We consider the $\FF_p$-algebra $(H^*(G,\FF_p),+\cup )$. We have the
isomorphism of graded $\FF_p$-algebras
$$H^*(G,\FF_p)\cong\begin{cases}S(V)&p=2\\
\lambda (V)\otimes S(V)&p>2\end{cases}.$$

We denote by $\tau$ the isomorphism $S(V)\to H^*(G,\FF_2)$, when $p=2$,
or $\Lambda (V)\otimes S(V)\to H^*(G,\FF_p)$ when $p>2$.

If $p=2$ then $V={\rm Hom}(G,\FF_2)$ identifies as $H^1(G,\FF_2)$ and
$\tau :S(V)\to H^*(G,\FF_2)$ is given by $f_1\cdots
f_n\mapsto f_1\cup\cdots\cup f_n\in H^n(G,\FF_2)$ $\forall
f_1,\ldots,f_n\in V$.

If $p>2$ then $V={\rm Hom}(G,\FF_p)$ from $\Lambda (V)$ identifies
with $H^1(G,\FF_p)$, while $V$ from $S(V)$ identifies with the image
of the Bockstein boundary map $\beta :H^1(G,\FF_p)\to H^2(G,\FF_p)$,
which is injective. So we may regard the isomorphism $\tau$ as being
defined as $\tau: \Lambda (V)\otimes S(\beta (V))\to H^*(G,\FF_p)$, by
$f_1\wedge\cdots\wedge f_l\otimes\beta (g_1)\cdots\beta (g_k)\mapsto
f_1\cup\cdots\cup f_l\cup\beta (g_1)\cup\cdots\cup\beta (g_k)\in
H^{2k+l}(G,\FF_p)$ $\forall f_1,\ldots f_l,g_1,\ldots,g_k\in V$. Note
that the homogeneous component of degree $n$ of $\Lambda (V)\otimes
S(\beta (V))$ is
$$(\Lambda (V)\otimes S(\beta (V)))^n=\bigoplus_{2k+l=n}\Lambda^l(V)\otimes
S^k(\beta (V)).$$

Equivlalently, this writes as
$$H^*(G,\FF_p)\cong\begin{cases}\FF_2[x_1,\ldots,x_r]&p=2\\
\lambda (x_1,\ldots,x_r)\otimes\FF_p[y_1,\ldots,y_r]&p>2\end{cases},$$
where $y_i=\beta (x_i)$. (See, e.g., Corollary II.4.3 and Theorem
II.4.4 in [AM].)

If $p=2$ then the isomorphism
$\FF_2[x_1,\ldots,x_r]\to H^*(G,\FF_2)$ is given by $x_{i_1}\cdots
x_{i_n}\mapsto x_{i_1}\cup\cdots\cup x_{i_n}\in H^n(G,\FF_2)$.

If $p>3$ then the isomorphism $\Lambda
(x_1,\ldots,x_r)\otimes\FF_p[y_1,\ldots,y_r]\to H^*(G,\FF_p)$ is given
by $x_{i_1}\wedge\cdots\wedge x_{i_l}\otimes y_{i_1}\cdots
y_{i_k}\mapsto x_{i_1}\cup\cdots\cup x_{i_l}\cup y_{i_1}\cup\cdots\cup
y_{i_k}\in H^{2k+l}(G,\FF_p)$.

In this paper we give explicit formulas for the reverse isomorphism,
$\tau^{-1}$. The elements of $H^*(G,\FF_p)$ will be given in terms of
normalized cochains, which allows the use of the $\II$-cochains
defined in \S 1.
\medskip

For $1\leq i\leq r$ we denote by $t_i=s_i-1$.

\blm (i) $t_1^{k_1}\cdots t_r^{k_r}$, with $0\leq k_i\leq p-1$, are a
basis of $\ZZ [G]$.

(ii) $t_1^{k_1}\cdots t_r^{k_r}$, with $0\leq k_i\leq p-1$ and
$(k_1,\ldots,k_r)\neq (0,\ldots,0)$ are a basis of $\II$.
\elm
\pf (i) If $k=(k_1,\ldots,k_r)$ with $k_i\in\ZZ_{\geq 0}$ then we denote
by $s^k=s_1^{k_1}\cdots s_r^{k_r}$ and by $t^k=t_1^{k_1}\cdots
t_r^{k_r}$. Since $s_1,\ldots,s_r$ are a basis of $G$ over $\FF_p$, we
have $G=\{ s^k\,\mid\, k\in S\}$, where $S=\{ 0,\ldots,p-1\}^r$. 

If $k=(k_1,\ldots,k_r)\in S$ then we define $\Sigma k=k_1+\cdots
+k_r$. We write the set $S$ as a sequence $k^1,\ldots,k^{p^r}$ such
that $\Sigma k^1\geq\cdots\geq\Sigma k^{p^r}$. Hence if $\Sigma
k<\Sigma l$ then $k=k^\alpha$ and $l=k^\beta$ for some $\alpha,\beta$
with $\alpha>\beta$.

If $k=(k_1,\ldots,k_r)\in S$ then
$t^k=\prod_{i=1}^rt_i^{k_i}=\prod_{i=1}^r(s_i-1)^{k_i}=s^k+$ a linear
combination of $s^l$ with $\Sigma l<\Sigma k$. So if $k=k^\alpha$ then
$t^{k^\alpha}$ writes as $s^{k^\alpha}+$ a linear combination $s^l$
with $\Sigma l<\Sigma k^\alpha$, which implies that $l\in\{ k^{\alpha
+1},\ldots, k^{p^r}\}$.

We consider the column matices $\underline
s=(s^{k^1},\ldots,s^{k^{p^r}})^T$ and $\underline
t=(t^{k^1},\ldots,t^{k^{p^r}})^T$. Since for every $\alpha$ we have
$t^{k^\alpha} =s^{k^\alpha}+$ a linear combination of $s^{k^\beta}$
with $\beta >\alpha$, we have $\underline t=A\underline s$ for some
$A\in M_{p^r}(\ZZ )$ which is upper unitriangular, so invertible. Then,
since $G=\{ s^k\,\mid\, k\in S\} =\{ s^{k^1},\ldots,s^{k^{p^r}}\}$ is a
basis for $\ZZ [G]$, this implies that $\{ t^k\,\mid\, k\in S\} =\{
t^{k^1},\ldots,t^{k^{p^r}}\}$ is a basis for $\ZZ [G]$.

(ii) From (i) we have that an element $\alpha\in\ZZ [G]$ writes
uniquely as $\alpha =\sum_{k\in S}c_kt^k$, with $c_k\in\ZZ$. Since
$t_i\in\II$ we have $t^k\in\II$, so $\varepsilon (t^k)=0$, for every
$k\in S\setminus\{ (0,\ldots,0)\}$. If $k=(0,\ldots,0)$ then
$\varepsilon (t^k)=\varepsilon (1)=1$. Hence $\varepsilon (\alpha
)=c_{0,\ldots,0}$. It follows that $\alpha\in\II =\ker\varepsilon$ iff
$c_{0,\ldots,0}=0$, i.e. if $\alpha$ has the form  $\alpha =\sum_{k\in
S\setminus\{ (0,\ldots,0)\}}c_kt^k$, with $c_k\in\FF_p$. Thus $\{
t^k\,\mid\, k\in S\setminus\{ (0,\ldots,0)\}\}$ is a basis for
$\II$. \qed

Recall that the Bockstein map $\beta :H^1(G,\FF_p)\to H^2(G,\FF_p)$ is
the connecting morphism resulting from the exact sequence
$0\to\FF_p\xrightarrow p\ZZ/p^2\ZZ\to\FF_p\to 0$. It coincides with
the reduction modulo $p$ of the connecting morphism $\tilde\beta
:H^1(G,\FF_p)\to H^2(G,\ZZ )$ resulting from the exact sequence
$0\to\ZZ\xrightarrow p\ZZ\to\FF_p\to 0$. 

\blm If $\tilde y_i=\tilde\beta (x_i)\in H^2(G,\ZZ )$, then $\tilde
y_i=[\tilde z_i]$, where $\tilde z_i\in C^n(G,\ZZ )$ is given by
$$\tilde z_i(u,v)=\left[\frac{k_i+l_i}p\right]
=\begin{cases}1&k_i+l_i\geq p\\
0&k_i+l_i<p\end{cases},$$
whenever $u=s_1^{k_1}\cdots s_r^{k_r}$ and $v=s_1^{l_1}\cdots
s_r^{l_r}$, with $0\leq k_1,\ldots,k_r,l_1,\ldots,l_r\leq p-1$.

Also $y_i=\beta (x_i)\in H^2(G,\FF_p)$ is given by $y_i=[z_i]$, where
$z_i\in C^n(G,\FF_p)$ is the reduction modulo $p$ of $\tilde z_i\in
C^n(G,\ZZ)$. 
\elm
\pf If $u=s_1^{k_1}\cdots s_r^{k_r}\in G$, with $0\leq k_h\leq p-1$,
then, since $x_i(s_j)=\delta_{i,j}$, we have $x_i(u)=k_i\in\FF_p$. In
the preimage of $x_i$ with respect to the map $C^1(G,\ZZ )\to
C^1(G,\FF_p)$ we have $\tilde x_i$, given by $\tilde x_i(u)=k_i\in\{
0,\ldots,p-1\}$.

If $u=s_1^{k_1}\cdots s_r^{k_r},v=s_1^{l_1}\cdots s_r^{l_r}\in G$,
with $0\leq u_h,v_h\leq p-1$, then $uv=s_1^{k_1+l_1}\cdots
s_r^{k_r+l_r}=s_1^{j_1}\cdots s_r^{j_r}$, where $0\leq j_h\leq p-1$,
$j_h=k_h+l_h-p[\frac{k_h+l_h}p]$. It follows that $d_1\tilde
x_i(u,v)=x_i(v)-x_i(uv)+x_i(u)=l_i-j_i+k_i=p[\frac{k_i+l_i}p]=p\tilde
z_i(u,v)$, so $d_1\tilde x_i=p\tilde z_i$. Then $\tilde z_i$ is in the
preimage of $d_1\tilde x_i$ with respect to the map $C^2(G,\ZZ
)\xrightarrow pC^2(G,\ZZ )$. It follows that $\tilde\beta
(x_i)=[\tilde z_i]$, as claimed. \qed

\bco If $p=2$ then $\beta (x)=x\cup x$ $\forall x\in H^1(G,\FF_2)$.
\eco
\pf If $x=0$ the statement is trivial. Otherwise we change the basis
$x_1,\ldots,x_r$ of $V$, together with the dual basis $s_1,\ldots,s_r$,
such that $x_1=x$. So we must prove that $y_1=\beta (x_1)=x_1\cup
x_1$. Let $u=s_1^{k_1}\cdots s_r^{k_r}$, $v=s_1^{l_1}\cdots s_r^{l_r}$
with $0\leq k_i,l_i\leq 1$. Then $x_1\cup x_1=[z'_1]$ where $z'_1\in
C^2(G,\FF_2)$ is given by $z'_1(u,v)=-x_1(u)x_1(v)=k_1l_1$. (We are in
characteristic $2$ so the minus sign can be ignored.) On the other
hand, $y_1=[z_1]$, with $z_1(u,v)=[\frac{k_1+l_1}2]$. Since
$k_1,l_1\in\{ 0,1\}$, we have $z_1(u,v)=1$ if $k_1=l_1=1$ and
$z_1(u,v)=0$ otherwise. Hence $z_1(u,v)=k_1l_1=z'_1(u,v)$. It follows
that $y_1=[z_1]=[z'_1]=x_1\cup x_1$. \qed

\blm If $f\in C_\II^n(G,\FF_p)$ and $\alpha,\beta\in T^n(\II )$ such
that $\alpha\equiv\beta\pmod p$, then $f(\alpha )=f(\beta )$.
\elm
\pf We have $\alpha -\beta =p\gamma$ for some $\gamma\in T^n(\II )$ so
$f(\alpha )-f(\beta )=pf(\gamma )=0$. \qed

\blm (i) If $1\leq i\leq r$ and $k\geq 1$, then $t^k=\sum_{h\leq
k}(-1)^{k-h}\binom kh(s_i^h-1)$.

(ii) In $\II$ we have $t_i^p\equiv 0\pmod p$ and
$t_i^{p-1}\equiv\sum_{h=1}^{p-1}(t^h-1)\pmod p$.

Consequently, if $f\in C_\II^n(G,\FF_p)$ and
$\alpha_1,\ldots,\alpha_n\in\II$ such that $\alpha_h=t_i^k$ for some
$1\leq h\leq n$, $1\leq i\leq r$ and $k\geq p$, then
$f(\alpha_1\otimes\cdots\otimes\alpha_n)=0$.
\elm
\pf (i) We have
$t_i^k=(s_i-1)^k=(s_i-1)^k-(1-1)^k=\sum_{h\leq k}(-1)^{k-h}\binom
khs_i^h-\sum_{h\leq k}(-1)^{k-h}\binom kh=\sum_{h\leq
k}(-1)^{k-h}\binom kh(s_i^h-1)$.

(ii) Since $p\mid\binom ph$ for $1\leq h\leq p-1$ we have
$t_i^p=\sum_{h=0}^p(-1)^{p-h}\binom ph(s_i^h-1)\equiv
(s_i^p-1)+(-1)^p(1-1)=0\pmod p$.

Over $\FF_p$ we have $\sum_{h\leq p-1}(-1)^{p-1-h}\binom{p-1}hX^h
=(X-1)^{p-1}=\frac{(X-1)^p}{X-1}=\frac{X^p-1}{X-1}=X^{p-1}+\cdots +1$
so $(-1)^{p-1-h}\binom{p-1}h\equiv 1\pmod p$ for $0\leq h\leq p-1$. It
follows that
$t_i^{p-1}=\sum_{h=0}^{p-1}(-1)^{p-1-h}\binom{p-1}h(s_i^h-1)
\equiv\sum_{h=0}^{p-1}(s_i^h-1)\pmod p$. The term for 
$h=0$ is zero so it can be ignored.

Since $t_i^p\equiv 0\pmod p$ we have $t_i^k\equiv 0\pmod p$ for $k\geq
p$. Hence if $\alpha_h=t_i^k$ for some $h$, then
$\alpha_1\otimes\cdots\otimes\alpha_n\equiv 0\pmod p$ and, by Lemma
2.4, $f(\alpha_1\otimes\cdots\otimes\alpha_n)=0$ for every $f\in
C_\II^n(G,\FF_p)$. \qed

\blm For every $1\leq i\leq r$ we denote by $f_i$, $\tilde g_i$ and
$g_i$ the correspondents of $x_i\in H^1(G,\FF_p)$, $\tilde y_i\in
H^1(G,\ZZ )$ and $y_i\in H^1(G,\FF_p)$ in $H_\II^1(G,\FF_p)$,
$H_\II^2(G,\ZZ )$ and $H_\II^2(G,\FF_p)$.

Let $\alpha,\beta$ in the basis of $\II$ from Lemma 2.1(ii).

$$(i)~ f_i(\alpha )=\begin{cases}1&\alpha =t_i\\
0&\text{otherwise}\end{cases}.$$

(ii) We have $\tilde g_i=[\tilde h_i]$, where $\tilde h_i\in
Z_\II^2(G,\ZZ )$ is given by
$$\tilde h_i(\alpha\otimes\beta
)=\begin{cases}\sum_{h=p}^{k+l}(-1)^{k+l-h}\binom{k+l}h&(\alpha,\beta
) =(t_i^k,t_i^l),~ k+l\geq p\\
0&\text{otherwise}\end{cases}.$$

(iii) We have $g_i=[h_i]$, where $h_i\in Z_\II^2(G,\ZZ )$ is given by
$$h_i(\alpha\otimes\beta
)=\begin{cases}1&(\alpha,\beta ) =(t_i^k,t_i^l),~ k+l=p\\
0&\text{otherwise}\end{cases}.$$
\elm
\pf (i) We have $f_i(u-1)=x_i(u)$ $\forall u\in G$. If $j\neq i$ then
for every $u\in G$ we have $x_i(s_ju)=x_i(s_j)+x_i(u)=x_i(u)$ and so
$f_i(t_ju)=f_i((s_j-1)u)=f_i(s_ju-1)-f_i(u-1)=x_i(s_ju)-x_i(u)=0$. Since
$\ZZ [G]$ is spanned by $G$, by linearity, we have $f_i(t_j\alpha
)=0$ $\forall\alpha\in\ZZ [G]$.

Let now $\alpha$ be an element in the basis of $\II$ from Lemma
2.1(ii). If $\alpha$ is not a power of $t_i$ then $\alpha =t_j\alpha'$
for some $j\neq i$ and $\alpha'\in\ZZ [G]$, which implies that
$f(\alpha )=0$. Assume now that $\alpha =t_i^k$ for some $1\leq k\leq
p-1$. By Lemma 2.5(i), $\alpha =\sum_{h=1}^k(-1)^{k-h}\binom
kh(s_i^h-1)$ so
$$f_i(\alpha )=\sum_{h=1}^k(-1)^{k-h}\binom khf_i(s_i^h-1)
=\sum_{h=1}^k(-1)^{k-h}\binom khx_i(s_i^h)
=\sum_{h=1}^k(-1)^{k-h}\binom khh.$$ 
Thus $f_i(\alpha )=P(1)$, where
$$P(X)=\sum_{h=1}^k(-1)^{k-h}\binom khhX^{h-1}
=\left(\sum_{h=0}^k(-1)^{k-h}\binom
khX^h\right)'=((X-1)^k)'=k(X-1)^{k-1}.$$
We have $P(1)=1$ if $k=1$ and $P(1)=0$ if $k>1$. Hence $f_i(\alpha
)=1$ if $\alpha =t_i$ and $f_i(\alpha )=0$ otherwise.

(ii) We have $\tilde g_i=[\tilde h_i]$, where $\tilde h_i\in
C_\II^2(G,\ZZ )$ is the $\II$-cochain correponding to $\tilde z_i\in
C^2(G,\ZZ )$ from Lemma 2.2. We have $\tilde h_i((u-1)\otimes
(v-1))=\tilde z_i(u,v)$.

If $u,v\in G$, $u=s_1^{k_1}\cdots s_r^{k_r}$,  $v=s_1^{l_1}\cdots
s_r^{l_r}$ then, by Lemma 2.2, $\tilde z_i(u,v)$ depends only on
$k_i$ and $l_i$. It follows that for every $j\neq i$ we have $\tilde
z_i(s_ju,v)=\tilde z_i(u,s_jv)=\tilde z_i(u,v)$ $\forall u,v\in G$.

If $u,v\in G$ and $j\neq i$ then $t_ju=(s_j-1)u=(s_ju-1)-(u-1)$ so
$t_ju\otimes (v-1)=(s_ju-1)\otimes (v-1)-(u-1)\otimes (v-1)$. It
follows that $\tilde h_i(t_ju\otimes (v-1))=\tilde h_i((s_ju-1)\otimes
(v-1))-\tilde h_i((u-1)\otimes (v-1))=\tilde z_i(s_ju,v)-\tilde
z_i(u,v)=0$. Since $u$, with $u\in G$, spans $\ZZ [G]$ and $v-1$, with
$v\in G$, spans $\II$, we have by linearity $\tilde
h_i(t_j\alpha\otimes\beta )=0$ $\forall\alpha\in\ZZ [G]$,
$\beta\in\II$. Similarly $\tilde h_i(\alpha\otimes t_j\beta )=0$
$\forall\alpha\in\II$, $\beta\in\ZZ [G]$.

Let now $\alpha,\beta$ be elements in the basis of $\II$ from Lemma
2.1(ii). If $\alpha$ is not a power of $t_i$ then $\alpha =t_j\alpha'$
so $\tilde h_i(\alpha\otimes\beta )=\tilde h_i(t_j\alpha'\otimes\beta
)=0$. Similarly when $\beta$ is not a power of $t_i$.

Assume now that $\alpha =t_i^k$, $\beta =t_i^l$, with $1\leq k,l\leq
p-1$. We have
$$\begin{aligned}
\tilde h_i(t_i^k\otimes t_i^l)&=\tilde
h_i\Big(\Big(\sum_{a\leq k}(-1)^{k-a}\binom
ka(s_i^a-1)\Big)\otimes\Big(\sum_{b\leq l}(-1)^{l-b}\binom
lb(s_i^b-1)\Big)\Big)\\
&=\sum_{a\leq k,b\leq l}(-1)^{k+l-a-b}\binom ka\binom lb\tilde
h_i((s_i^a-1)\otimes (s_i^b-1))\\
&=\sum_{a\leq k,b\leq l}(-1)^{k+l-a-b}\binom ka\binom lb\tilde
z_i(s_i^a\otimes s_i^b)\\
&=\sum_{a\leq k,b\leq l,a+b\geq p}(-1)^{k+l-a-b}\binom ka\binom
lb.
\end{aligned}$$
(Here we used the formula $\tilde z_i(s_i^a\otimes
s_i^b)=[\frac{a+b}p]$ from Lemma 2.2.)

If $k+l<p$ then the sum above is empty so $\tilde h_i(t_i^k\otimes
t_i^l)=0$. If $k+l\geq p$ then $h:=a+b$ in the sum above takes
values between $p$ and $k+l$. If we use the convention that $\binom
ak=0$ for every integer $a$ outside the interval $[0,k]$ and similarly
for $\binom bl$ then
$$\tilde h_i(t_i^k\otimes
t_i^l)=\sum_{h=p}^{k+l}\sum_{a+b=h}(-1)^{k+l-h}\binom ka\binom
lb=\sum_{h=p}^{k+l}(-1)^{k+l-h}\binom{k+l}h.$$
This concludes the proof of (ii).

(iii) We prove that $h_i$ is the reduction modulo $p$ of $\tilde
h_i$. We only have to prove that $h_i(\alpha\otimes\beta )$ is the
reduction modulo $p$ of $\tilde h_i(\alpha\otimes\beta )$ in the case
when $(\alpha,\beta )=(t_i^k,t_i^l)$, with $k+l\geq p$. Since $k,l\leq
p-1$ we have $k+l=p+c$ with $0\leq c\leq p-2$. If we make the
substitution $h\to h+p$, we get
$$\tilde h_i(t_i^k\otimes t_i^l)=\sum_{h=p}^{p+c}(-1)^{p+c-h}\binom{p+c}h=
\sum_{h=0}^c(-1)^{c-h}\binom{p+c}{p+h}.$$
But, by Lucas's theorem, since $0\leq h,c\leq p-1$ we have
$\binom{p+c}{p+h}\equiv\binom ch\binom 11=\binom ch\pmod p$. Hence
$\tilde h_i(t_i^k\otimes t_i^l)\equiv\sum_{h=0}^c(-1)^{c-h}\binom
ch\pmod p$, which is $1$ if $c=0$, i.e. if $k+l=p$, and it is $0$
otherwise. Hence the conclusion. \qed

To solve our problem, we first state it in terms of $\II$-cocycles. We
denote by $\tau_\II$ the isomorphism with values in $H_\II^*(G,\FF_p)$
induced by $\tau$ via the isomorphism $H^*(G,\FF_p)\cong
H_\II^*(G,\FF_p)$. Since the corresondents of $x_i,y_i\in
H^*(G,\FF_p)$ in $H_\II^*(G,\FF_p)$ are $f_i,g_i$, in the case $p=2$,
the isomorphism $\tau_\II :\FF_2[x_1,\ldots,x_r]\to H_\II^*(G,\FF_2)$
is given by $x_{i_1}\cdots x_{i_n}\mapsto f_{i_1}\cup\cdots\cup
f_{i_n}$ and, in the case $p>2$, $\tau_\II :\Lambda
(x_1,\ldots,x_r)\otimes\FF_p[y_1,\ldots,y_r]\to H_\II^*(G,\FF_p)$
is given by $x_{i_1}\wedge\cdots\wedge x_{i_l}\otimes y_{j_1}\cdots
y_{j_k}\mapsto f_{i_1}\cup\cdots\cup f_{i_l}\cup g_{j_1}\cup\cdots\cup
g_{j_k}$. To find $\tau^{-1}$ we first determine $\tau_\II^{-1}$. 

We now prove the main result in the case $p=2$. Since we are in
characteristic $2$, the powers of $-1$ can be ignored.

\btm If $p=2$ then
$\tau_\II^{-1}:H_\II^*(G,\FF_2)\to\FF_2[x_1,\ldots,x_r]$ is given by
$$[f]\mapsto\sum_{1\leq i_1,\ldots,i_n\leq
r}f(t_{i_1}\otimes\cdots\otimes t_{i_n})x_{i_1}\cdots x_{i_n}$$
for every $f\in Z_\II^n(G,\FF_2)$.
\etm
\pf First we define $\theta
:Z_\II^*(G,\FF_2)\to\FF_2[x_1,\cdots,x_r]$, given for $f\in
Z_\II^n(G,\FF_2)$ by $f\mapsto\sum_{1\leq i_1,\ldots,i_n\leq
r}f(t_{i_1}\otimes\cdots\otimes t_{i_n})x_{i_1}\cdots x_{i_n}$ and we
prove that $\theta (B_\II^*(G,\FF_2))=0$.

Let $g\in C_\II^{n-1}(G,\FF_2)$. We use Remark 1.3 and ignore the
powers of $-1$. Then
$$\begin{aligned}
\theta (d_{n-1}g) &
=\sum_{i_1,\ldots,i_n}d_{n-1}g(t_{i_1}\otimes\cdots\otimes
t_{i_n})x_{i_1}\cdots x_{i_n}\\
{} & =\sum_{i_1,\ldots,i_n}\bigg(\sum_{h=1}^{n-1}
g(t_{i_1}\otimes\cdots\otimes t_{i_h}t_{i_{h+1}}\otimes\cdots\otimes
t_{i_n}\bigg) x_{i_1}\cdots x_{i_n}.
\end{aligned}$$

Note that if $i_h=i_{h+1}$ then $t_{i_h}t_{i_{h+1}}=t_{i_h}^2$ so, by
Lemma 2.5(ii), $g(t_{i_1}\otimes\cdots\otimes
t_{i_h}t_{i_{h+1}}\otimes\cdots\otimes t_{i_n})=0$. Therefore 
$$\theta (d_{n-1}g)
=\sum_{i_1,\ldots,i_n}\sum_{h=1}^{n-1}{}\raisebox{3pt}'g(t_{i_1}\otimes\cdots\otimes
t_{i_h}t_{i_{h+1}}\otimes\cdots\otimes t_{i_n})x_{i_1}\cdots
x_{i_n},$$
where $\sum_{h=1}'^{n-1}$ is the sum restricted to indices $h$ with
$i_h\neq i_{h+1}$.

Note that the sum above contains terms where we have factors of the
type $f(t_{j_1},\otimes\cdots\otimes
t_{j_l}t_{j_{l+1}}\otimes\cdots\otimes t_{j_n})$, with
$j_l<j_{l+1}$. Each factor of this type has the form appears in
exactly two terms, corresponding to $k=l$ and
$(i_1,\ldots,i_n)=(j_1,\ldots,j_n)$ or
$(j_1,\ldots,j_{l+1},j_l,\ldots,j_n)$. But both these terms are equal
to $f(t_{j_1},\otimes\cdots\otimes
t_{j_l}t_{j_{l+1}}\otimes\cdots\otimes t_{j_n})x_{j_1}\cdots x_{j_n}$ so
they cancel each other. Therefore $\theta (d_{n-1}g)=0$ and we have
$\theta (B_\II^n(G,\FF_2))=0$. More generally $\theta
(B_\II^*(G,\FF_2))=0$. It follows that $\theta$ induces a
morphism $\bar\theta :H_\II^*(G,\FF_2)\to\FF_2[x_1,\ldots,x_r]$, given
by $\bar\theta ([f])=\theta (f)$ $\forall f\in Z_\II^n(G,\FF_2)$.

We now prove that $\bar\theta\tau_\II$ is the identity on
$\FF_2[x_1,\ldots,x_r]$, so that $\tau_\II^{-1}=\bar\theta$. We prove
that $\bar\theta\tau_\II (\eta )=\eta$ for every monomial $\eta
=x_{j_1}\cdots x_{j_n}$. By Lemma 1.4(ii), with the power of $-1$
ignored, $\tau_\II (\eta )=f_{j_1}\cup\cdots\cup f_{j_n}=[f]$, where
$f\in C_\II^n(G,\FF_2)$,
$f(\alpha_1,\otimes\cdots\otimes\alpha_n)=f_{j_1}(\alpha_1)\cdots
f_{j_n}(\alpha_n)$ $\forall\alpha_1,\ldots,\alpha_n\in\II$. Then
$\bar\theta\tau_\II (\eta )=\bar\theta
([f])=\sum_{i_1,\ldots,i_n}f(t_{i_1}\otimes\cdots\otimes
t_{i_n})x_{i_1}\cdots x_{i_n}$. But $f(t_{i_1}\otimes\cdots\otimes
t_{i_n})=f_{j_1}(t_{i_1})\cdots f_{j_n}(t_{i_n})$, which, by Lemma
2.6(i), is equal to $1$ if $i_h=j_h$ $\forall h$ and it is equal to
$0$ otherwise. Hence $\bar\theta\tau (\eta )=x_{j_1}\cdots
x_{j_n}=\eta$. \qed

\bco If $p>2$ then $\tau^{-1}:H^*(G,\FF_2)\to\FF_2[x_1,\ldots,x_n]$ is
given by
$$[a]\mapsto\sum_{1\leq i_1,\ldots,i_n\leq
r}a(s_{i_1},\ldots,s_{i_n})x_{i_1}\cdots x_{i_n}$$
for every $a\in Z^n(G,\FF_2)$.
\eco

\pf We have $\tau^{-1}([a])=\tau_\II^{-1}([f])$, where $f\in
C_\II^n(G,\FF_2)$ is the $\II$-cocycle corresponding to $a\in
C^n(G,\FF_2)$. Then
$a(s_{i_1},\ldots,s_{i_n})=f((s_{i_1}-1)\otimes\cdots\otimes
(s_{i_n}-1))=f(t_{i_1}\otimes\cdots\otimes t_{i_n})$ for every $1\leq
i_1,\ldots,i_n\leq r$. So, by using the formula for
$\tau_\II^{-1}([f])$ from Theorem 2.7, we get the desired result. \qed

The case $p>2$ is more complicated we need some preliminaries.

For any $x\in V$ and $l\geq 0$ we denote by $x^{\wedge
l}\in\Lambda^l(V)$ the wedge product of $l$ copies of $x$. We have
$x^{\wedge 0}=1\in\FF_p=\Lambda^0(V)$, $x^{\wedge 1}=x\in
V=\Lambda^1(V)$ and for $l\geq 2$ we have $x^{\wedge
  l}=x\wedge\cdots\wedge x=0$.

If $x\in V$ and $m\geq 0$ the we define $x^{(m)}\in (\Lambda
(V)\otimes S(\beta (V)))^m$ by $x^{(m)}=x^{\wedge l}\otimes\beta
(x)^k\in\Lambda^l(V)\otimes S^k(\beta (V))$, where $m=2k+l$, with
$l\in\{ 0,1\}$, i.e. $x^{(m)}=1\otimes\beta (x)^k$ if $m=2k$ and
$x^{(m)}=x\otimes\beta (x)^k$ if $m=2k+1$. It is easy to see that
$x^{(m)}x^{(m')}=x^{(m+m')}$ if $mm'$ is even, but $x^{(m)}x^{(m')}=0$
if $mm'$ is odd.

In the particular case when $x=x_i$ we have $\beta (x_i)$ so
$$x_i^{(m)}=\begin{cases}1\otimes y_i^k&m=2k\\
x_i\otimes y_i^k&m=2k+1\end{cases}.$$

Note that a basis of $\Lambda (V)\otimes S(\beta (V))=\Lambda
(x_1,\ldots,x_r)\otimes\FF_p[y_1,\ldots,y_r]$ is made of products
$\eta =x_{i_1}\wedge\cdots\wedge x_{i_l}\otimes y_1^{k_1}\cdots
y_r^{k_r}$, with $0\leq l\leq r$, $1\leq i_1<\cdots <i_l\leq r$ and
$k_i\geq 0$ $\forall i$. Now the wedge product
$x_{i_1}\wedge\cdots\wedge x_{i_l}$ writes as $x_1^{\wedge
l_1}\wedge\cdots\wedge x_r^{\wedge l_r}$, where $l_i=1$ if $i\in\{
i_1,\ldots,i_l\}$ and $l_i=0$ otherwise. Then we have
$$\eta =x_1^{\wedge l_1}\wedge\cdots\wedge x_r^{\wedge l_r}\otimes
y_1^{k_1}\cdots y_r^{k_r}=(x_1^{\wedge l_1}\otimes y_1^{k_1})\cdots
(x_r^{\wedge l_r}\otimes y_r^{k_r}) =x_1^{(n_1)}\cdots x_r^{(n_r)},$$
where $n_i=2k_i+l_i$. Moreover, since $x_i^{(n_i)}$ is homogeneous of
degree $n_i$, $\eta$ is homogeneous of degree $n_1+\cdots +n_r$. Hence
we have:

\blm The set $\{ x_1^{(n_1)}\cdots x_r^{(n_r)}\,\mid\, n_1+\cdots
+n_r=n\}$ is a basis for\\ $(\Lambda
(x_1,\ldots,x_r)\otimes\FF_p[y_1,\ldots,y_r])^n$.
\elm

\bdf For every $1\leq i\leq r$ and every $m\geq 0$ w define
$t_{i,m}\in T^m(\II )$ by
$$t_{i,m}=\begin{cases}(t_i^{p-1}\otimes t_i)^{\otimes k}& m=2k\\
t_i\otimes (t_i^{p-1}\otimes t_i)^{\otimes k}& m=2k+1\end{cases}.$$
\edf

\blm If $1\leq i\leq r$ and $m\geq 0$ then the $\tau_\II
(x_i^{(m)})=[f_{i,m}]$, with $f_{i,m}\in C_\II^m(G,\FF_p)$ such that
for every $\alpha_1,\ldots,\alpha_m$ in the basis of $\II$ from Lemma
2.1(ii) we have

(i) $f_{i,m}(\alpha_1\otimes\cdots\otimes\alpha_m)=1$ if 
$(\alpha_1,\ldots,\alpha_m)=(t_i^{q_1},\ldots,t_i^{q_m})$, with $1\leq
q_j\leq p-1$ $\forall j$, such that either $m$ is even and
$q_1+q_2=q_3+q_4=\cdots =q_{m-1}+q_m=p$ or $m$ is odd, $q_1=1$ and
$q_2+q_3=q_4+q_5=\cdots =q_{m-1}+q_m=p$.

In particular, $f_{i,m}(t_{i,m})=1$.

(ii) $f_{i,m}(\alpha_1\otimes\cdots\otimes\alpha_m)=0$ otherwise.
\elm
\pf If $m=2k$ then $\tau_\II (x_i^{(m)})=\tau_\II (1\otimes
y_i^k)=g_i\cup\cdots\cup g_i$, while if $m=2k+1$ then $\tau_\II
(x_i^{(m)})=\tau_\II (x\otimes y_i^k)=f_i\cup g_i\cup\cdots\cup
g_i$. Since $f_i=[f_i]$ and $g_i=[h_i]$, by Lemma 1.4(ii), we have
$\tau_\II (x_i^{(m)})=[f_{i,m}]$, with $f_{i,m}\in C_\II^m(G,\FF_p)$
such that for every $\alpha_1,\ldots,\alpha_m\in\II$
we have
$$f_{i,m}(\alpha_1\otimes\cdots\otimes\alpha_m)=\begin
{cases}h_i(\alpha_1\otimes\alpha_2)h_i(\alpha_3\otimes\alpha_4)\cdots
h_i(\alpha_{m-1}\otimes\alpha_m)&m\text{ is }even\\
f_i(\alpha_1)h_i(\alpha_2\otimes\alpha_3)h_i(\alpha_4\otimes\alpha_5)\cdots
h_i(\alpha_{m-1}\otimes\alpha_m)&m\text{ is }odd\end{cases}.$$
(Note that the $(-1)^{\frac{l(l-1)}2}$ factor from Lemma 1.4(ii) is in
both cases $1$ since the corresponding $n_1,\ldots,n_r$ sequence is
$2,\ldots,2$ or $1,2,\ldots,2$, so $l=0$ or $1$, respectively.)

By Lemma 2.6, if $\alpha_1,\ldots,\alpha_m$ belong to the basis of
$\II$ from Lemma 2.1(ii), then all factors of
$f_{i,m}(\alpha_1\otimes\ldots\otimes\alpha_m)$ are either $0$ or $1$
so $f_{i,m}(\alpha_1\otimes\ldots\otimes\alpha_m)=0$ or $1$. More
precisely, we have $h_i(\alpha_j\otimes\alpha_{j+1})=1$ iff
$(\alpha_j,\alpha_{j+1})=(t_i^{q_j},t_i^{q_{j+1}})$, where $1\leq
q_j,q_{j+1}\leq p-1$ and $q_j+q_{j+1}=p$ and we have
$f_i(\alpha_1)=1$ iff $\alpha_1=t_i$. It follows that
$f_{i,m}(\alpha_1\otimes\ldots\otimes\alpha_m)=1$ precisely in the
cases described in (i) and it equals $0$ otherwise.

Finally, note that $t_{i,m}=\alpha_1\otimes\cdots\otimes\alpha_m$,
with $(\alpha_1,\ldots,\alpha_m)=(t_i^{q_1},\ldots,t_i^{q_m})$, where
$(q_1,\ldots,q_m)=(p-1,1,p-1,1,\ldots,p-1,1)$ if $m$ is even and
$(q_1,\ldots,q_m)=(1,p-1,1,p-1,1,\ldots,p-1,1)$ if $m$ is odd. Hence
$\alpha_1,\ldots,\alpha_m$ satisfy the conditions of (i) and we have
$f_{i,m}(t_{i,m})=1$. \qed

We consider the action of the symmetric group $S_n$ on $C^n(G,\FF_p)$
given by
$$\sigma
a(u_1,\ldots,u_n)=\sgn (\sigma
)a(u_{\sigma^{-1}(1)},\ldots,u_{\sigma^{-1}(n)})~\forall u_1,\ldots,u_n\in G.$$

If $f\in C_\II^n(G,\FF_p)$ is the $\II$-cochain corresponding to $a\in
C^n(G,\FF_p)$ then for every $\sigma\in S_n$ the $\II$-cochain
corresponding to $\sigma a$ is
$\sigma f$, where 
$$\sigma f(\alpha_1\otimes\cdots\otimes\alpha_n)=\sgn (\sigma
)f(\alpha_{\sigma^{-1}(1)}\otimes\cdots\otimes\alpha_{\sigma^{-1}(n)})
~\forall\alpha_1,\ldots,\alpha_n\in\II.$$

Next, if $n_1+\cdots +n_r=n$ then we denote by $Sh(n_1,\ldots,n_r)$
the set of all $(n_1,\ldots,n_r)$-shuffles, i.e.
$$Sh(n_1,\ldots,n_r)=\{\sigma\in S_n\,\mid\,\sigma (h)<\sigma
(h+1)\,\forall h,\, h\neq n_1+\cdots +n_i\,\forall 1\leq i<r\}.$$ 
Equivalently, the condition from the definition of
$Sh(n_1,\ldots,n_r)$ can be written as
$$\sigma (n_1+\cdots +n_{i-1}+1)<\cdots <\sigma (n_1+\cdots
+n_i)~\forall 1\leq i\leq r.$$

\btm If $p>2$ then $\tau_\II^{-1}:H_\II^*(G,\FF_p)\to\Lambda
(x_1,\ldots,x_r)\otimes\FF_p[y_1,\ldots,y_r]$ is given by
$$[f]\mapsto\sum_{n_1+\cdots +n_r=n}
c_{n_1,\ldots,n_r}x_1^{(n_1)}\cdots x_r^{(n_r)},$$
for every $f\in Z_\II^n(G,\FF_p)$, where 
$$c_{n_1,\ldots,n_r}=(-1)^{\frac{l(l-1)}2}\sum_{\sigma\in
Sh(n_1,\ldots,n_r)}\sigma f(t_{1,n_1}\otimes\cdots\otimes
t_{r,n_r}),$$
with $l=|\{ i\,\mid\, 1\leq i\leq r,\, n_i\text{ is odd}\}|$.
\etm
\pf Same as for Theorem 2.7, we first define the function $\theta
:Z_\II^*(G,\FF_p)\to\Lambda
(x_1,\ldots,x_r)\otimes\FF_p[y_1,\ldots,y_r]$ by
$f\mapsto\sum_{n_1+\cdots +n_r=n}c_{n_1,\ldots,n_r}x_1^{(n_1)}\cdots
x_r^{(n_r)}$ $\forall f\in Z_\II^n(G,\FF_p)$. We prove that $\theta
(B_\II^*(G,\FF_p))=0$, i.e. that $\theta (f)=0$ if $f=d_{n-1}g$ for
some $g\in C_\II^{n-1}(G,\FF_p)$. We must prove that each coefficient
$c_{n_1,\ldots,n_r}$ is zero.

Let $\alpha_1\otimes\cdots\otimes\alpha_n
=t_{1,n_1}\otimes\cdots\otimes t_{r,n_r}$. Since $t_{i,n_i}\in
T^{n_i}(\II )$ $\forall i$, we have
$$\alpha_{n_1+\cdots +n_{i-1}+1}\otimes\cdots\otimes\alpha_{n_1+\cdots
+n_i}=t_{i,n_i}.$$
Hence the sequence $\alpha_{n_1+\cdots
+n_{i-1}+1},\ldots,\alpha_{n_1+\cdots +n_i}$ is
$t_i^{p-1},t_i,\ldots,t_i^{p-1},t_i$ or\\
$t_i,t_i^{p-1},t_i,\ldots,t_i^{p-1},t_i$, depending on the parity of
$n_i$. Then we have
$$\begin{aligned}
c_{n_1,\ldots,n_r}&=(-1)^{\frac{l(l-1)}2}\sum_{\sigma\in
Sh(n_1,\ldots,n_r)}\sigma f(\alpha_1\otimes\cdots\otimes\alpha_n)\\
&=(-1)^{\frac{l(l-1)}2}\sum_{\sigma\in  Sh(n_1,\ldots,n_r)}\sgn
(\sigma)f(\alpha_{\sigma^{-1}(1)}\otimes\cdots\otimes\alpha_{\sigma^{-1}(n)}).
\end{aligned}$$
Since $f=d_{n-1}g$, we have $f(\beta_1\otimes\cdots\otimes\beta_n)
=\sum_{h=1}^{n-1}(-1)^hg(\beta_1\otimes\cdots\otimes\beta_h
\beta_{h+1}\otimes\cdots\otimes\beta_r)$ for every
$\beta_1,\ldots,\beta_n\in\II$. It follows that
$c_{n_1,\ldots,n_r}=(-1)^{\frac{l(l-1)}2}\sum_{h=1}^{n-1}(-1)^h\Sigma_h$,
where
$$\Sigma_h=\sum_{\sigma\in Sh(n_1,\ldots,n_r)}\sgn (\sigma
)g(\alpha_{\sigma^{-1}(1)}\otimes\cdots\otimes\alpha_{\sigma^{-1}(h)}
\alpha_{\sigma^{-1}(h+1)}\otimes\cdots\otimes\alpha_{\sigma^{-1}(r)}).$$
We prove that $\Sigma_h=0$ $\forall h$.

For convenience, for every $\sigma\in S_n$ we denote by
$$\psi (\sigma )=\sgn (\sigma
)g(\alpha_{\sigma^{-1}(1)}\otimes\cdots\otimes\alpha_{\sigma^{-1}(h)}
\alpha_{\sigma^{-1}(h+1)}\otimes\cdots\otimes\alpha_{\sigma^{-1}(r)}).$$
Then we have $\Sigma_h=\sum_{\sigma\in Sh(n_1,\ldots,n_r)}\psi
(\sigma )$. We denote by $A$ the set of all $\sigma\in
Sh(n_1,\ldots,n_r)$ such that both $\sigma^{-1}(h)$ and
$\sigma^{-1}(h+1)$ belong to the same interval $[n_1+\cdots
+n_{i-1}+1,n_1+\cdots +n_i]$ for some $1\leq i\leq r$ and we denote
by $B=Sh(n_1,\ldots,n_r)\setminus A$. Then
$\Sigma_h=\Sigma'_h+\Sigma''_h$, where $\Sigma'_h=\sum_{\sigma\in 
A}\psi (\sigma )$ and $\Sigma''_h=\sum_{\sigma\in B}\psi (\sigma
)$. We prove that $\Sigma'_h=\Sigma''_h=0$.

We denote by $h'=\sigma^{-1}(h)$ and $h''=\sigma^{-1}(h+1)$, so that
$\sigma (h')=h$ and $\sigma (h'')=h+1$.

If $\sigma\in A$ then $h',h''\in [n_1+\cdots +n_{i-1}+1,n_1+\cdots
+n_i]$ for some $i$. Since $\sigma (n_1+\cdots +n_{i-1}+1)<\cdots
<\sigma (n_1+\cdots +n_i)$ and $\sigma (h'')=\sigma (h')+1$, we must
have $h''=h'+1$. Since the sequence $\alpha_{n_1+\cdots
+n_{i-1}+1},\ldots,\alpha_{n_1+\cdots +n_i}$ is
$t_i^{p-1},t_i,\ldots,t_i^{p-1},t_i$ or
$t_i,t_i^{p-1},t_i,\ldots,t_i^{p-1},t_i$, we have
$(\alpha_{h'},\alpha_{h''})=(\alpha_{h'},\alpha_{h'+1})=(t_i^{p-1},t_i)$
or $(t_i,t_i^{p-1})$. In both cases,
$\alpha_{\sigma^{-1}(h)}\alpha_{\sigma^{-1}(h+1)}
=\alpha_{h'}\alpha_{h''}=t_i^p$. By Lemma 2.5(ii), it follows that
$g(\alpha_{\sigma^{-1}(1)}\otimes\cdots\otimes\alpha_{\sigma^{-1}(h)}
\alpha_{\sigma^{-1}(h+1)}\otimes\cdots\otimes\alpha_{\sigma^{-1}(r)})=0$. Hence
$\psi (\sigma )=0$ $\forall\sigma\in A$, which implies $\Sigma'_h=0$.

If $\sigma\in B$ then $h'\in [n_1+\cdots +n_{i'-1}+1,n_1+\cdots
+n_{i'}]$ and $h''\in [n_1+\cdots +n_{i''-1}+1,n_1+\cdots +n_{i''}]$
for some $i'\neq i''$.

We claim that $\tau_h\sigma\in B$, where $\tau_h\in S_n$ is the
transposition $(h,h+1)$. First we prove that $\tau_h\sigma\in
Sh(n_1,\ldots,n_r)$. Since $\tau_h$ only permutates $h$ and $h+1$, we
have $\tau_h\sigma (j)=\sigma (j)$ for all $j$ except
$j=\sigma^{-1}(h)=h'$ and $j=\sigma^{-1}(h+1)=h''$. We have
$\tau_h\sigma (h')=\tau_h(h)=h+1=\sigma (h')+1$,
$\tau_h\sigma (h'')=\tau_h(h+1)=h=\sigma (h'')-1$ and
$\tau_h\sigma (j)=\sigma (j)$ for $j\neq h',h''$. It follows that for
$i\neq i',i''$ the inequalities $\sigma (n_1+\cdots +n_{i-1}+1)<\cdots
<\sigma (n_1+\cdots +n_i)$ from the definition of $Sh(n_1,\ldots,n_r)$
remain the same if we replace $\sigma$ by $\tau_h\sigma$. If $i=i'$ or
$i''$ they change to
$$\sigma (n_1+\cdots +n_{i'-1}+1)<\cdots <\sigma (h')+1<\cdots
<\sigma (n_1+\cdots +n_{i'}),$$
$$\sigma (n_1+\cdots +n_{i''-1}+1)<\cdots <\sigma (h'')-1<\cdots
<\sigma (n_1+\cdots +n_{i''}).$$
The only inequalities that may fail to hold are 
$\sigma (h')+1<\sigma (h'+1)$ (if $h'<n_1+\cdots +n_{i'}$)
and $\sigma (h''-1)<\sigma (h'')-1$ (if $n_1+\cdots
+n_{i''-1}+1<h''$). We have $\sigma (h')<\sigma (h'+1)$ so if
$\sigma (h')+1<\sigma (h'+1)$ fails then $\sigma (h'+1)=\sigma
(h')+1=h+1=\sigma (h'')$ so $h''=h'+1\in [n_1+\cdots
+n_{i'-1}+1,n_1+\cdots +n_{i'}]$. Contradiction. Similarly, $\sigma
(h''-1)<\sigma (h'')$ so if $\sigma (h''-1)<\sigma (h'')-1$ fails then
$\sigma (h''-1)=\sigma (h'')-1=h=\sigma (h')$ so $h'=h''-1\in
[n_1+\cdots +n_{i''-1}+1,n_1+\cdots +n_{i''}]$. Again,
contradiction. Thus $\tau_h\sigma\in Sh(n_1,\ldots,n_r)$. We have
$\tau_h\sigma (h')=h+1$ and $\tau_h\sigma (h'')=h$ so $(\tau_h\sigma
)^{-1}(h)=h''\in [n_1+\cdots +n_{i''-1}+1,n_1+\cdots +n_{i''}]$ and 
$(\tau_h\sigma )^{-1}(h+1)=h'\in [n_1+\cdots +n_{i'-1}+1,n_1+\cdots
+n_{i'}]$. Since $i''\neq i'$ we have $\tau_h\sigma\in B$.

As a consequence, $B$ writes as a disjoint union of right cosets
$\langle\tau_h\rangle\sigma =\{\sigma,\tau_h\sigma\}$,
$B=\bigsqcup_{\sigma\in C}\{\sigma,\tau_h\sigma\}$. Thus
$\Sigma''_h=\sum_{\sigma\in C}(\psi (\sigma )+\psi (\tau_h\sigma
))$. For every $\sigma\in C$ we have
$(\tau_h\sigma)^{-1}=\sigma^{-1}\tau_h$ so
$(\tau_h\sigma)^{-1}(h)=\sigma^{-1}(h+1)$,
$(\tau_h\sigma)^{-1}(h+1)=\sigma^{-1}(h)$ and if $j\neq h,h+1$ then
$(\tau_h\sigma^{-1})(j)=\sigma^{-1}(j)$. We also have
 $\sgn (\tau_h\sigma )=-\sgn (\sigma )$. It follows that
\begin{multline*}
\psi (\tau_h\sigma )=\sgn (\tau_h\sigma
)g(\alpha_{(\tau_h\sigma
)^{-1}(1)}\otimes\cdots\otimes\alpha_{(\tau_h\sigma )^{-1}(h)}
\alpha_{(\tau_h\sigma
)^{-1}(h+1)}\otimes\cdots\otimes\alpha_{(\tau_h\sigma )^{-1}(r)})\\
=-\sgn (\sigma
)g(\alpha_{\sigma^{-1}(1)}\otimes\cdots\otimes\alpha_{\sigma^{-1}(h+1)}
\alpha_{\sigma^{-1}(h)}\otimes\cdots\otimes\alpha_{\sigma^{-1}(r)})=-\psi
(\sigma ).
\end{multline*}
Thus $\Sigma''_h=\sum_{\sigma\in C}(\psi (\sigma )+\psi (\tau_h\sigma
))=\sum_{h\in C}0=0$.

Since $\theta (B_\II^*(G,\FF_p))=0$, we a morphism $\bar\theta
:H_\II^*(G,\FF_p)\to\Lambda
(x_1,\ldots,x_r)\otimes\FF_p[y_1,\ldots,y_r]$, given by $\bar\theta
([f])=\theta (f)$ $\forall f\in Z_\II^n(G,\FF_p)$. We prove that
$\bar\theta\tau_\II$ is the identity of $\Lambda
(x_1,\ldots,x_r)\otimes\FF_p[y_1,\ldots,y_r]$, so that
$\tau_\II^{-1}=\bar\theta$. We prove that $\bar\theta\tau_\II (\eta
)=\eta$ for the monomials $\eta =x_1^{(m_1)}\cdots x_r^{(m_r)}$, which
generate $\Lambda (x_1,\ldots,x_r)\otimes\FF_p[y_1,\ldots,y_r]$. Let
$m_1+\cdots +m_r=n$. By Lemma 2.10,
$\tau_\II (x_i^{(m_i)})=[f_{i,m_i}]\in H_\II^{m_i}(G,\FF_p)$ and so
$\tau_\II (\eta )=[f]$, where $f\in Z_\II^n(G,\FF_p)$,
$f=f_{1,m_1}\cup\cdots\cup f_{r,m_r}$. Since $f_{i,m_i}\in
Z_\II^{m_i}(G,\FF_p)$, we have by Lemma 1.4(ii) 
$$f(\beta_1\otimes\cdots\otimes\beta_n)
=(-1)^{\frac{l(l-1)}2}\prod_{i=1}^rf_{i,m_i}(\beta_{m_1+\cdots
+m_{i-1}+1}\otimes\cdots\otimes\beta_{m_1+\cdots +m_i}),$$
where $l=|\{ i\,\mid\, 1\leq i\leq r,\, m_i\text{ is odd}\}|$.

For every $n_1,\ldots,n_r$ with $n_1+\cdots +n_r=n$ we denote by
$c_{n_1,\ldots,n_r}$ the coefficient of $x_1^{(n_1)}\cdots
x_r^{(n_r)}$ in $\bar\theta\tau_\II (\eta)=\bar\theta ([f])$. Same as
above, we consider $\alpha_1,\ldots,\alpha_n\in\II$ such that
$\alpha_1\otimes\cdots\otimes\alpha_n=t_{1,n_1}\otimes\cdots\otimes
t_{r,n_r}$. For every $1\leq i\leq n_i$ we have that
$\alpha_{n_1+\cdots +n_{i-1}+1},\ldots,\alpha_{n_1+\cdots +n_i}$ are
powers of $t_i$ such that $\alpha_{n_1+\cdots
+n_{i-1}+1}\otimes\cdots\otimes\alpha_{n_1+\cdots +n_i}=t_{i,n_i}$. So
the sequence $\alpha_1,\ldots,\alpha_n$ contains $n_i$ powers of $t_i$
for every $i$. Same happens for the sequnece $\alpha_{\sigma
(1)},\ldots,\alpha_{\sigma (n)}$ for every $\sigma\in S_n$. We have
$$\begin{aligned}
c_{n_1,\ldots,n_r}&=(-1)^{\frac{l'(l'-1)}2}\sum_{\sigma\in
Sh(n_1,\ldots,n_r)}\sigma f(\alpha_1\otimes\cdots\otimes\alpha_n)\\
&=(-1)^{\frac{l'(l'-1)}2}\sum_{\sigma\in
Sh(n_1,\ldots,n_r)}\sgn (\sigma
)f(\alpha_{\sigma^{-1}(1)}\otimes\cdots\otimes\alpha_{\sigma (n)})\\
&=(-1)^{\frac{l'(l'-1)}2}\sum_{\sigma\in
Sh(n_1,\ldots,n_r)}\sgn (\sigma )(-1)^{\frac{l(l-1)}2}\cdot\\
&\qquad\qquad\qquad\qquad\cdot\prod_{i=1}^r
f_{i,m_i}(\alpha_{\sigma^{-1}(m_1+\cdots +m_{i-1}+1)}
\otimes\cdots\otimes \alpha_{\sigma^{-1}(m_1+\cdots +m_i)}),
\end{aligned}$$
where $l'=|\{ i\,\mid\, 1\leq i\leq r,\, n_i\text{ is odd}\}|$.

Recall that every $\alpha_j$ has the form $t_i$ or $t_i^{p-1}$ for
some $i$, so it belongs to the basis of $\II$ from Lemma 2.1(ii). By
Lemma 2.10, it follows that $f_{i,m_i}(\alpha_{\sigma^{-1}(m_1+\cdots +m_{i-1}+1)}
\otimes\cdots\otimes \alpha_{\sigma^{-1}(m_1+\cdots +m_i)})=0$ unless
$\alpha_{\sigma^{-1}(m_1+\cdots
+m_{i-1}+1)},\ldots,\alpha_{\sigma^{-1}(m_1+\cdots +m_i)}$ are powers
of $t_i$. Hence if
$\prod_{i=1}^rf_{i,m_i}(\alpha_{\sigma^{-1}(m_1+\cdots +m_{i-1}+1)}
\otimes\cdots\otimes\alpha_{\sigma^{-1}(m_1+\cdots +m_i)})\neq 0$
then $\forall i$ the sequence $\alpha_{\sigma
(1)},\ldots\alpha_{\sigma (n)}$ contains $n_i$ powers of $t_i$, which
implies that $n_i=m_i$ $\forall i$. Consequently,
$c_{n_1,\ldots,n_r}=0$ if $(n_1,\ldots,n_r)\neq (m_1,\ldots,m_r)$.

If $(n_1,\ldots,n_r)=(m_1,\ldots,m_r)$ then $l'=l$ so the factors
$(-1)^{\frac{l(l-1)}2}$ and $(-1)^{\frac{l'(l'-1)}2}$ cancel each
other. We get
$$c_{m_1,\ldots,m_r}=\sum_{\sigma\in Sh(m_1,\ldots,m_r)}\sgn
(\sigma )\prod_{i=1}^rf_{i,m_i}(\alpha_{\sigma^{-1}(m_1+\cdots +m_{i-1}+1)}
\otimes\cdots\otimes \alpha_{\sigma^{-1}(m_1+\cdots +m_i)}).$$
As seen above, if the term corresponding to $\sigma$ in the sum above is
$\neq 0$ then for every $i$ we have that
$\alpha_{\sigma^{-1}(m_1+\cdots
+m_{i-1}+1)},\ldots,\alpha_{\sigma^{-1}(m_1+\cdots +m_i)}$ are powers
of $t_i$. But the only powers of $t_i$ in the sequence
$\alpha_1,\ldots,\alpha_n$ are $\alpha_{m_1+\cdots
+m_{i-1}+1},\ldots,\alpha_{m_1+\cdots +m_i}$. It follows that
$\{\sigma^{-1}(m_1+\cdots +m_{i-1}+1),\ldots,\sigma^{-1}(m_1+\cdots
+m_i)\} =\{ m_1+\cdots +m_{i-1}+1,\ldots,m_1+\cdots +m_i\}$ and so $\{
m_1+\cdots +m_{i-1}+1,\ldots,m_1+\cdots +m_i\} =\{\sigma (m_1+\cdots
+m_{i-1}+1),\ldots,\sigma (m_1+\cdots +m_i)\}$. But $\sigma\in
Sh(m_1,\ldots,m_r)$ so $\sigma (m_1+\cdots +m_{i-1}+1)<\cdots <\sigma
(m_1+\cdots +m_i)$. It follows that $\sigma (j)=j$ for $m_1+\cdots
+m_{i-1}+1\leq j\leq m_1+\cdots +m_i$ and for every $1\leq i\leq
r$. Hence $\sigma =1$. So the only non-zero term of
$c_{m_1,\ldots,m_r}$ corresponds to $\sigma =1$ and, by Lemma 2.10(i),
we have
$$c_{m_1,\ldots,m_r}=\prod_{i=1}^rf_{i,m_i}(\alpha_{m_1+\cdots
+m_{i-1}+1}\otimes\cdots\otimes\alpha_{m_1+\cdots
+m_i})=\prod_{i=1}^rf_{i,m_i}(t_{i,m_i})=1.$$
In coclusion, $\bar\theta\tau_\II (\eta )=\sum_{n_1+\cdots
+n_r=n}c_{n_1,\ldots,n_r}x_1^{(n_1)}\cdots
x_1^{(n_1)}=x_1^{(m_1)}\cdots x_1^{(m_1)}=\eta$. \qed

\bdf If $1\leq i\leq r$, $m\geq 0$, $k=[m/2]$ and $q_1,\ldots,q_k$ are
nonnegative integers then we define
$s_{i,m,q_1,\ldots,q_k}\in G^m$ by
$$s_{i,m,q_1,\ldots,q_k}=\begin{cases}(s_i^{q_1},s_i,\ldots
s_i^{q_k},s_i)& m=2k\\
(s_i,s_i^{q_1},s_i,\ldots s_i^{q_k},s_i)& m=2k+1\end{cases}$$
and $t_{i,m,q_1,\ldots,q_k}\in T^m(\II )$ by
$t_{i,m,q_1,\ldots,q_k}=(u_1-1)\otimes\cdots (u_m-1)$, where
$(u_1,\ldots,u_m)=s_{i,m,q_1,\ldots,q_k}$, i.e.
$$t_{i,m,q_1,\ldots,q_k}=\begin{cases}(s_i^{q_1}-1)\otimes
(s_i-1)\otimes\cdots\otimes (s_i^{q_k}-1)\otimes (s_i-1))& m=2k\\
(s_i-1)\otimes (s_i^{q_1}-1)\otimes (s_i-1)\otimes\cdots\otimes
(s_i^{q_k}-1)\otimes (s_i-1))& m=2k+1\end{cases}.$$
\edf

\blm If $m\geq 0$ is an integer and $k=[m/2]$ then
$$t_{i,m}\equiv\sum_{1\leq q_1,\ldots,q_k\leq
p-1}t_{i,m,q_1,\ldots,q_k}\pmod p.$$
\elm
\pf We have $t_i=s_i-1$ and, by Lemma 2.5(ii),
$t_i^{p-1}\equiv\sum_{q=1}^{p-1}(s_i^q-1)\pmod p$. So if $m=2k$ then
$$\begin{aligned}
t_{i,m}&=t_i\otimes t_i^{p-1}\otimes\cdots\otimes t_i\otimes t_i^{p-1}\\
&\equiv
(s_i-1)\otimes\sum_{q_1=1}^{p-1}(s_i^{q_1}-1)\otimes\cdots\otimes
(s_i-1)\otimes\sum_{q_k=1}^{p-1}(s_i^{q_k}-1)\\
&=\sum_{1\leq q_1,\ldots,q_k\leq p-1}(s_i-1)\otimes
(s_i^{q_1}-1)\otimes\cdots\otimes (s_i-1)\otimes (s_i^{q_k}-1)\\
&=\sum_{1\leq q_1,\ldots,q_k\leq p-1}t_{i,m,q_1,\ldots,q_k}\pmod p.
\end{aligned}$$
The case $m=2k+1$ follows from the case $m=2k$ by noting that
$t_{i,2k+1}=t_i\otimes t_{i,2k}=(s_i-1)\otimes t_{i,2k}$ and
$t_{i,2k+1,q_1,\ldots,q_k}=(s_i-1)\otimes
t_{i,2k,q_1,\ldots,q_k}$. \qed

\bco If $p>2$ then $\tau^{-1}:\Lambda
(x_1,\ldots,x_r)\otimes\FF_p[y_1,\ldots,y_r]\to H^*(G,\FF_p)$ is given
by
$$[a]\mapsto\sum_{n_1+\cdots
+n_r=n}c_{n_1,\ldots,n_r}x_1^{(n_1)}\cdots x_r^{(n_r)}$$
for every $a\in Z^n(G,\FF_p)$, where
$$\begin{aligned}
c_{n_1,\ldots,n_r}=(-1)^{\frac{l(l-1)}2}\sum_{\sigma\in
Sh(n_1,\ldots,n_r)}\sum_{1\leq q_{i,j}\leq p-1}\sigma
a(s_{1,n_1,q_{1,1},\ldots,q_{1,k_1}},\ldots,s_{r,n_r,q_{r,1},\ldots,q_{r,k_r}}),
\end{aligned}$$
with $l=|\{ i\,\mid\, 1\leq i\leq r,\, n_i\text{ is odd}\}|$ and
$k_i=[n_i/2]$.

Here by the sum $\sum_{1\leq q_{i,j}\leq p-1}$ we mean that every
variable $q_{i,j}$, with $1\leq i\leq r$ and $1\leq j\leq k_i$, takes
values between $1$ and $p-1$.

Also $s_{1,n_1,q_{1,1},\ldots,q_{1,k_1}},\ldots,s_{r,n_r,q_{r,1},\ldots,q_{r,k_r}}$
is the concatenation of the sequences
$s_{i,m_i,q_{i,1},\ldots,q_{i,k_i}}$ for $1\leq i\leq r$, of lengths
$n_1,\ldots,n_r$.
\eco
\pf Let $f\in Z_\II^n(G,\FF_p)$ be the $\II$-cochain corresponding to
$a\in Z^n(G,\FF_p)$. Then for every $\sigma\in S_n$ the $\II$-cochain
corresponding to $\sigma a$ is $\sigma f$. Hence\\ $\sigma
f((u_1-1)\otimes\cdots\otimes (u_n-1))=\sigma a(u_1,\ldots,u_n)$
$\forall u_1,\ldots,u_n\in G$.

We have $\tau^{-1}([a])=\tau_\II^{-1}([f])$. To prove our statement,
we write the terms of $\tau_\II^{-1}([f])$ in terms of $a$.

By Lemma 2.12, $t_{1,n_1}\otimes\cdots\otimes t_{r,n_r}$ is congruent
modulo $p$ to
\begin{multline*}
\Big(\sum_{1\leq q_{1,1},\ldots,q_{1,k_1}\leq p-1}
t_{1,n_1,q_{1,1},\ldots,q_{1,k_1}}\Big)\otimes\cdots\otimes\Big(\sum_{1\leq
q_{r,1},\ldots,q_{r,k_r}\leq p-1} t_{r,n_r,q_{r,1},\ldots,q_{r,k_r}}\Big)\\
=\sum_{1\leq q_{i,j}\leq
p-1}t_{1,n_1,q_{1,1},\ldots,q_{1,k_1}}\otimes\cdots\otimes
t_{r,n_r,q_{r,1},\ldots,q_{r,k_r}}
\end{multline*}
By Lemma 2.4, this implies that
$$\sigma f(t_{1,n_1}\otimes\cdots\otimes t_{r,n_r})=\sum_{1\leq
q_{i,j}\leq p-1}\sigma
f(t_{1,n_1,q_{1,1},\ldots,q_{1,k_1}}\otimes\cdots\otimes
t_{r,n_r,q_{r,1},\ldots,q_{r,k_r}}).$$
By definition, if $s_{i,m,q_1,\ldots,q_k}=(u_1,\ldots,u_m)$ then
$t_{i,m,q_1,\ldots,q_k}=(u_1-1)\otimes\cdots\otimes (u_m-1)$. It
follws that if
$(s_{1,n_1,q_{1,1},\ldots,q_{1,k_1}},\ldots,s_{r,n_r,q_{r,1},\ldots,q_{r,k_r}})
=(u_1,\ldots,u_n)$ then $t_{1,n_1,q_{1,1},\ldots,q_{1,k_1}}\otimes\cdots\otimes
t_{r,n_r,q_{r,1},\ldots,q_{r,k_r}}=(u_1-1)\otimes\cdots\otimes
(u_n-1)$. This implies that
$\sigma
a(s_{1,n_1,q_{1,1},\ldots,q_{1,k_1}},\ldots,s_{r,n_r,q_{r,1},\ldots,q_{r,k_r}})
=\sigma f(t_{1,n_1,q_{1,1},\ldots,q_{1,k_1}}\otimes\cdots\otimes
t_{r,n_r,q_{r,1},\ldots,q_{r,k_r}})$ $\forall\sigma\in S_n$. Therefore
$$\sigma f(t_{1,n_1}\otimes\cdots\otimes t_{r,n_r})=\sum_{1\leq
q_{i,j}\leq p-1}\sigma a(s_{1,n_1,q_{1,1},\ldots,q_{1,k_1}},\ldots,
s_{r,n_r,q_{r,1},\ldots,q_{r,k_r}})$$
so our result follows from Theorem 2.11. \qed

{\bf Remarks}

1. In the case $p=2$, the formulas for $\tau_\II^{-1}([f])$ and
$\tau^{-1}([a])$ from Theorem 2.7 and Corollary 2.8 can be written in
terms of $ x_1^{n_1}\cdots x_r^{n_r}$, with $n_1+\cdots +n_r=n$, which
are a basis of $\FF_p[x_1,\ldots,x_r]^n$. As a consequence of Theorem
2.7, we have
$$\tau_\II^{-1}([f])=\sum_{n_1+\cdots
+n_r=n}c_{n_1,\ldots,n_r}x_1^{n_1}\cdots x_r^{n_r},$$
with
$$c_{n_1,\ldots,n_r}=\sum_{(i_1,\ldots,i_n)\in
S(n_1,\ldots,n_r)}f(t_{i_1}\otimes\cdots\otimes t_{i_n}),$$
where $S(n_1,\ldots,n_r)=\{ (i_1,\ldots,i_n)\,\mid x_{i_1}\cdots
x_{i_n}=x_1^{n_1}\cdots x_r^{n_r}\}$.

We now extend Definition 2 to the case $p=2$. Since $p-1=1$, we have
$t_{i,m}=t_i^{\otimes m}\in T^m(\II )$, regardless of the parity of
$m$. We prove that $c_{n_1,\ldots,n_r}$ is given by the same formula
from Theorem 2.11.

We have $(i_1,\ldots,i_n)\in S(n_1,\ldots,n_r)$ iff for every $1\leq
i\leq r$ the sequence $i_1,\ldots,i_n$ contains $n_i$ copies of
$i$. Let $j_{i,1}<\cdots <j_{i,n_i}$ be the $n_i$ indices $j$ such
that $i_j=i$. Then we have a bijection $\psi :S(n_1,\ldots,n_r)\to
Sh(n_1,\ldots,n_r)$ given by $(i_1,\ldots,i_n)\mapsto\sigma$, where
$\sigma$ is defined on each interval $[n_1+\cdots
+n_{i-1}+1,n_1+\cdots, +n_i]$ by\\ $\sigma (n_1+\cdots
+n_{i-1}+h)=j_{i,h}$ $\forall 1\leq h\leq n_i$.

We have $t_{1,n_1}\otimes\cdots\otimes
t_{r,n_r}=\alpha_1\otimes\cdots\otimes\alpha_n$, where for each $i$ we
have $\alpha_{n_1+\cdots +n_{i-1}+1}=\cdots =\alpha_{n_1+\cdots
+n_i}=t_i$ so that $\alpha_{n_1+\cdots
+n_{i-1}+1}\otimes\cdots\otimes\alpha_{n_1+\cdots +n_i}=t_i^{\otimes
n_i}=t_{i,n_i}$. Let $(i_1,\ldots,i_n)\in S(n_1,\ldots,n_r)$ and let
$\sigma =\psi (i_1,\ldots,i_n)$, as above. Let $1\leq j\leq
n$. We have $j=j_{i,h}$ for some $1\leq i\leq r$ and $1\leq h\leq
n_i$. Then $i_j=i$ and $\sigma (n_1+\cdots +n_{i-1}+h)=j_{i,h}=j$ so
$\sigma^{-1}(j)=n_1+\cdots +n_{i-1}+h\in [n_1+\cdots
+n_{i-1}+1,n_1+\cdots +n_i]$, which implies that
$\alpha_{\sigma^{-1}(j)}=t_i=t_{i_j}$. It follows that
$f(t_{i_1}\otimes\cdots\otimes t_{i_n})
=f(\alpha_{\sigma^{-1}(1)}\otimes\cdots\otimes\alpha_{\sigma^{-1}(n)})
=\sigma f(\alpha_1\otimes\cdots\otimes\alpha_n)=\sigma
f(t_{1,n_1}\otimes\cdots\otimes t_{r,n_r})$. (We are in
characteristic $2$ so the factor $\sgn (\sigma )$ from the definition
of $\sigma f$ can be ignored.) Since $\psi :S(n_1,\ldots,n_r)\to
Sh(n_1,\ldots,n_r)$ is a bijection, the formula for
$c_{n_1,\ldots,n_r}$ als writes as
$$c_{n_1,\ldots,n_r}=\sum_{\sigma\in Sh(n_1,\ldots,n_r)}\sigma
f(t_{1,n_1}\otimes\cdots\otimes t_{r,n_r}),$$
which coincides with the formula from the case $p>2$.

For $1\leq i\leq r$ and $m\geq 0$ we denote by $s_{i,m}\in G^m$,
$s_{i,m}=(s_i,\ldots,s_i)$. Recall that $t_{i,m}\in T^m(\II )$ is
given by $t_{i,m}=t_i\otimes\cdots\otimes
t_i=(s_i-1)\otimes\cdots\otimes (s_i-1)$. Hence if 
$(s_{1,n_1},\ldots,s_{r,n_r})=(u_1,\ldots,u_n)\in G^n$ then
$t_{1,n_1}\otimes\cdots\otimes t_{r,n_r}=(u_1-1)\otimes\cdots\otimes
(u_n-1)\in T^n(\II )$. Then if $a\in C^n(G,\FF_2)$ and $f\in
C_\II^n(G,\FF_2)$ is the corresponding $\II$-cochain, then for every
$\sigma\in S$ the $\II$-cochain correponding to $\sigma a$ is $\sigma
f$ and so $\sigma f(t_{1,n_1}\otimes\cdots\otimes t_{r,n_r})=\sigma
a(s_{1,n_1},\ldots,s_{r,n_r})$. Thus if $a\in Z^n(G,\FF_2)$ then the
coefficients $c_{n_1,\ldots,n_r}$ from
$\tau^{-1}([a])=\tau_\II^{-1}([f])$ write as
$$c_{n_1,\ldots,n_r}=\sum_{\sigma\in Sh(n_1,\ldots,n_r)}\sigma
a(s_{1,n_1},\ldots,s_{r,n_r}).$$
This is  the same formula from Corollary 2.13, since $s_{i,m}$ coincides
with $s_{i,m,1,\ldots,1}$ of Definition 3 and the sum $\sum_{1\leq
q_{i,j}\leq p-1}$ from Corollary 2.13 in the case $p=2$ has only one
term, with $q_{i,j}=1$ $\forall i,j$.

This is explained by the fact that the formula for $\tau (x^m)$ in the
case $p=2$ coincides to the formula for $\tau (x^{(m)})$ in the case
$p>2$. Indeed, if denote by $x^{\cup m}$ the cup product of $m$ copies
of $x$, then if $p=2$ we have $\tau (x^m)=x^{\cup m}$. When $p>2$ if
$m=2k$ then $\tau (x^{(m)})=\tau (1\otimes\beta (x)^k)=\beta (x)^{\cup
k}$, while if $m=2k+1$ then $\tau (x^{(m)})=\tau (x\otimes
y^k)=x\cup\beta (x)^{\cup k}$. But when $p=2$, by Corollary 2.3, we
have $\beta (x)=x\cup x$. So if $m=2k$ then $\beta (x)^{\cup k}=(x\cup
x)^{\cup k}=x^{\cup 2k}=x^{\cup m}=x^{\cup m}=\tau (x^m)$ and if
$m=2k+1$ then $x\cup\beta (x)^{\cup k}=x\cup (x\cup x)^{\cup
k}=x^{\cup 2k+1}=x^{\cup m}=x^{\cup m}=\tau (x^m)$.
\medskip

2. We now determine the number $N_n$ of terms
$a_{\alpha_1,\ldots,\alpha_n}$ involved in the formula for
$\tau^{-1}([a])$ from Corollary 2.13. For the term $c_{n_1,\ldots,n_r}$
we have $|Sh(n_1,\ldots,n_r)|=\binom n{n_1,\ldots,n_r}$ and the number
of choiches $1\leq q_{i,j}\leq p-1$ for $1\leq i\leq r$ and $1\leq
j\leq k_i=[n_i/2]$ is $(p-1)^{[n_1/2]+\cdots +[n_r/2]}$. So
$c_{n_1,\ldots,n_r}$ is the sum of $\binom
n{n_1,\ldots,n_r}(p-1)^{[n_1/2]+\cdots +[n_r/2]}$ terms. It follows
that
$$N_n=\sum_{n_1+\cdots +n_r=n}\binom
n{n_1,\ldots,n_r}(p-1)^{[n_1/2]+\cdots +[n_r/2]}.$$
Since $\frac 1{n!}\binom n{n_1,\ldots,n_r}=\frac 1{n_1!\cdots n_r!}$,
the formal series $F(X)=\sum_{n\geq 0}\frac{N_n}{n!}X^n$ writes as
$$\begin{aligned}
F(X)&=\sum_{n\geq 0}\sum_{n_1+\cdots +n_r=n}\frac 1{n_1!\cdots
n_r!}(p-1)^{[n_1/2]+\cdots +[n_r/2]}X^n\\
&=\sum_{n_1,\ldots,n_r}\frac 1{n_1!\cdots
n_r!}(p-1)^{[n_1/2]+\cdots +[n_r/2]}X^{n_1+\cdots +n_r}=G(X)^r,
\end{aligned}$$
where $G(X)=\sum_{n\geq 0}\frac 1{n!}(p-1)^{[n/2]}X^n$. 

We have $(p-1)^{[n/2]-n/2}=1$ if $n$ is even and $=\frac
1{\sqrt{p-1}}$ if $n$ is odd. So if $A_p=\frac 12(1+\frac
1{\sqrt{p-1}})$ and $B_p=\frac 12(1-\frac 1{\sqrt{p-1}})$, such that
$A_p+B_p=1$ and $A_p-B_p=\frac 1{\sqrt{p-1}}$, then
$(p-1)^{[n/2]-n/2}=A_p+(-1)^nB_p$. Hence
$(p-1)^{[n/2]}=\sqrt{p-1}^n(p-1)^{[n/2]-n/2}
=A_p\sqrt{p-1}^n+B_p(-\sqrt{p-1})^n$ and so
$$G(X)=\sum_{n\geq 0}\frac
1{n!}(A_p\sqrt{p-1}^n+B_p(-\sqrt{p-1})^n)
=A_pe^{\sqrt{p-1}X}+B_pe^{-\sqrt{p-1}X}.$$
It follows that $F(X)=G(X)^r=\sum_{k=0}^r\binom
rkA_p^{r-k}B_p^ke^{(r-2k)\sqrt{p-1}X}$.
By dentifying the coefficient of $X^n$ we get
$\frac{N_n}{n!}=\sum_{k=0}^r\binom
rkA_p^{r-k}B_p^k\frac{((r-2k)\sqrt{p-1})^n}{n!}$. Hence
$$N_n=\sum_{k=0}^r\binom rkA_p^{r-k}B_p^k((r-2k)\sqrt{p-1})^n$$

If $0\leq k\leq r$ then $|(r-2k)\sqrt{p-1}|\leq\sqrt{p-1}$, with
equality when $k=0$ or $r$. It follows that, as $n\to\infty$,
$$N_n\cong A_p^r(r\sqrt{p-1})^n+B_p^r(-r\sqrt{p-1})^n
=(A_p^r+(-1)^nB_p^r)(r\sqrt{p-1})^n.$$
By comparisson, the total number of coefficients $a_{u_1,\ldots,u_n}$,
with $u_1,\ldots,u_n\in G\setminus\{ 0\}$, which determine a cochain
$a\in C^n(G,\FF_p)$, is $(|G|-1)^n=(p^r-1)^n$.
\bigskip

{\bf References}
\bigskip

[AM] Adem, A. and Milgram, R.J., {\em Cohomology of Finite Groups},
Springer--Verlag Grundlehren {\bf 309} (2004).

[HS] Hilton, P.J. and Stammbach, U., {\em A Course in Homological
Algebra}, Graduate Texts in Mathematics, 2nd edition,
Springer--Verlag, New York 1997.
\bigskip

Institute of Mathematics Simion Stoilow of the Romanian
Academy,\\ Calea Grivitei 21, RO-010702 Bucharest, Romania.

E-mail address: Constantin.Beli@imar.ro

\end{document}